\documentclass[preprint,3p,12pt]{elsarticle}

\usepackage{lineno,hyperref}
\modulolinenumbers[5]

\usepackage{hyperref}

\usepackage{epstopdf}

\usepackage{mathrsfs}
\usepackage{amssymb}
\usepackage{amsfonts}
\usepackage{latexsym}
\usepackage{amsmath}
\usepackage{xcolor}
\usepackage{wrapfig}
\usepackage{floatflt}

\usepackage{graphicx}
\usepackage{subcaption}
\def\lb{\label}

\newcommand{\er}[1]{\textrm{(\ref{#1})}}

\newtheorem{theorem}{\bf Theorem}[section]

\def\a{\alpha}         
\def\b{\beta}          
\def\g{\gamma}         
  \def\cD{{\mathcal D}}       
\def\d{\delta}

    \def\cJ{{\mathcal J}}

\def\l{\lambda}

  \def\cR{{\mathcal R}}       
         
    \def\cT{{\mathcal T}}

\def\o{\omega}

       \def\vp{\varphi}    

\def\Z{{\mathbb Z}}    \def\R{{\mathbb R}}   \def\C{{\mathbb C}}

\def\lt{\biggl}                  \def\rt{\biggr}
               \def\wt{\widetilde}

\let\ge\geqslant                 \let\le\leqslant

\def\iy{\infty}
\def\sm{\setminus}               \def\es{\emptyset}
\def\ss{\subset}

\def\el2{\ell^{\,2}}             \def\1{1\!\!1}

\newtheorem{corollary}[theorem]{\bf Corollary}

\let\ge\geqslant
\let\le\leqslant

\newcommand{\ca}{\begin{cases}}
\newcommand{\ac}{\end{cases}}
\newcommand{\ma}{\begin{pmatrix}}
\newcommand{\am}{\end{pmatrix}}
\def\eq{\begin{equation}}
\def\qe{\end{equation}}
\renewcommand{\[}{\begin{equation}}
\renewcommand{\]}{\end{equation}}

\begin{document}

\begin{frontmatter}

\title{On some explicit integrals related to ``fractal foothills''}
\date{\today}

\author
{Anton A. Kutsenko}

\address{Jacobs University (International University Bremen), 28759 Bremen, Germany; email: akucenko@gmail.com}

\begin{abstract}
In previous papers, we attempted to analyze the complete loop counting functions that count all loops in an infinite random walk, represented by the digits of a real number. In this paper, the consideration will be restricted to the partial loop counting functions $V$ that count the returns to the origin only. This simplification allows us to find closed-form expressions for various integrals related to $V$. Some applications to the complete loop counting functions, in particular, their connections with Bernoulli polynomials, are also provided. 

\end{abstract}



\end{frontmatter}


{\section{Lead paragraph}\lb{sec0}}

Random walks, including those on graphs, are widely used in many branches of science. Of particular interest are precise analytical results of non-trivial characteristics, such as the weighted number of loops in a random walk. This is due to both the well-known difficult questions of the distribution of self-avoiding random walks and the statistical analysis of weather and climate changes. Weather change graphs are considered, for example, by Smilkov and Kocarev; they also present some analytical results. Random walks are conveniently represented by real variables, in which the digits correspond to the steps of the walk. The loop counting functions depend on these real variables. They have a fractal structure and contain the most complete information about the cycles in random walks. The more various analytical results can be obtained regarding these functions, the more accurately one can analyze the deep patterns of the distribution of cycles and their correlations in a random walk. In this article, we will obtain an exact expression for all sorts of integrals of these functions, including the Fourier transform. The integrals are expressed through the determinants of special Hessenberg matrices, through continued fractions, and through Bernoulli polynomials. Thus, even in the one-dimensional case, it is possible to obtain beautiful results that link different sections of mathematics. The operator approach presented in some integrals can be generalized to the multidimensional case.

{\section{Introduction}\lb{sec1}}
Any real number in its dyadic representation can be considered as an infinite random walk, where the digits correspond to the steps of the walk. One of the most interesting problems is the distribution of self-avoiding random walks, the walks that have no loops. Thus, let us provide a brief scheme of motivations starting from more simple and going to more complex objects: functions that count the number of returns to zero (fractal foothills) $\to$ loop counting functions (LCF) (fractal mountains) $\to$ self-avoiding random walks (SAW) as zeros of LCF $\to$ possible applications to various hard problems on distributions of SAW in a multidimensional case, see details in \cite{K}. But, of course, the main motivation should be the search for interesting relations between the objects, including various formulas based on continued fractions, determinants of special matrices, classical polynomials, etc. We will focus on the connections between these beautiful components of classical analysis and the stochastic curves mentioned above: "fractal foothills" and "fractal mountains."

Some useful information about random walks itself, including open questions, and representations of real numbers as random walks, is available in \cite{ABBB} and \cite{S}. The current work is a further development of some results presented in \cite{K}. In particular, unexpected relations with Bernoulli polynomials and determinants of Hessenberg matrices are found. While the work is motivated by \cite{K}, it can be read completely independently. All the results have an independent, complete form, understandable without any motivations and references to other literature. Let us start with the main results and postpone the further discussion to the end of the Introduction section. 

There is an interesting intersection between this topic and random walks on graphs. Namely, when the number of states is large, but the transition, due to physical reasons, is possible only between neighboring states, and we want to estimate the number of weighted loops (returns to some state, weighted depending on time), our function $V$ well approximates the corresponding measure. One such graph, representing the weather dynamics, is considered in \cite{SK}. The number of nodes of such a graph can be arbitrarily large, depending on the details of the information we need. One of the achievements mentioned in \cite{SK} is the ability to obtain analytical results related to random processes, such as random walks on graphs. It is worth noting that for some complex nonlinear process characteristics, analytical results can lead to some very beautiful mathematics.

Any $x\in[-1,1]$ except a countable set of some dyadic rationals can be uniquely expanded as
\[\lb{001}
 x=\frac{x_0}2+\frac{x_1}{2^2}+\frac{x_2}{2^3}+...,\ \ x_n\in\{-1,+1\}.
\]
For $\l\in\C$, $|\l|<1$, let us define the function that counts the number of returns to the origin multiplied by the exponential weight 
\[\lb{002}
 V(x)=1+\sum_{n=0}^{+\iy}\l^{n+1}L_n(x),\ \ \ L_n(x)=\ca 1,& \sum_{j=0}^nx_j=0,\\
                                                  0,& otherwise. \ac
\]
This function can be uniformly approximated by piecewise constant functions that are linear combinations of characteristic functions of intervals with dyadic endpoints. The function $V$ is even, measurable, and has a typical fractal structure, see Fig. \ref{fig1}. The function satisfies infinite number of symmetry relations: if $\wt x$ is $x$ with some swapped digits $x_{2n}\leftrightarrow x_{2n+1}$, see \er{001}, then $V(x)=V(\wt x)$. (It is important that $x_{2n}\leftrightarrow x_{2n+1}$, not $x_{2n+1}\leftrightarrow x_{2n+2}$.) 
\begin{figure}
  \centering
    \includegraphics[width=0.99\textwidth]{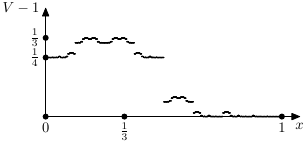}
  \caption{The plot of $V-1$ for $\l=1/2$. \lb{fig1}}  
\end{figure}

\begin{figure}
    \centering
    \begin{subfigure}[b]{0.49\textwidth}
        \includegraphics[width=\textwidth]{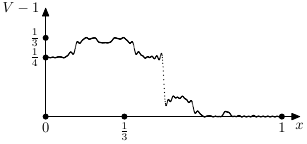}
        \caption{$100$ harmonics}
    \end{subfigure}  
    \begin{subfigure}[b]{0.49\textwidth}
        \includegraphics[width=\textwidth]{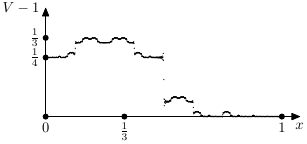}
        \caption{$500$ harmonics}
    \end{subfigure} 
    \caption{For $\l=1/2$, the approximation of $V-1$ in the trigonometric Fourier basis $\{\cos n\pi x\}_{n\ge0}$ is plotted.}\label{fig2}
\end{figure}

Let us assume by definition that $\sum_{i=a}^b\cdot=0$ and $\prod_{i=a}^b\cdot=1$ if $b<a$. The notation $|{\bf A}|$ for square matrices ${\bf A}$ means the determinant of ${\bf A}$. The binomial coefficients are denoted by $\binom{n}{m}$.  We formulate our main result.
\begin{theorem}\lb{T1} i) Let $P(x)=\sum_{n=0}^Np_nx^n$ be a polynomial with $p_n\in\C$. Then
\begin{multline}\lb{003}
 \int_{-1}^1P(V(x))dx=\\
 \frac2{\prod_{n=1}^N\sqrt{1-\l^{2n}}}
 \begin{vmatrix}
 1  & -\sqrt{1-\l^2} & 0 & ... & 0 & 0 \\
 1  & \frac{\binom{2}{1}\l^2}{1+\sqrt{1-\l^2}} & -\sqrt{1-\l^4} & ... & 0 & 0 \\
 1  & \frac{\binom{3}{1}\l^2}{1+\sqrt{1-\l^2}} & \frac{\binom{3}{2}\l^4}{1+\sqrt{1-\l^4}} & ... & 0 & 0 \\
 ... & ... & ... & ... & ... & ... \\
 1  & \frac{\binom{N}{1}\l^2}{1+\sqrt{1-\l^2}} & \frac{\binom{N}{2}\l^4}{1+\sqrt{1-\l^4}}  & ... & \frac{\binom{N}{N-1}\l^{2N-2}}{1+\sqrt{1-\l^{2N-2}}} & -\sqrt{1-\l^{2N}} \\
 p_0 & p_1 & p_2 & ... & p_{N-1} & p_N \end{vmatrix}.
\end{multline}
In particular 
\[\lb{003a}
 \int_{-1}^1V(x)^Ndx=
 \frac2{\prod_{n=1}^N\sqrt{1-\l^{2n}}}
 \begin{vmatrix}
 1  & -\sqrt{1-\l^2} & 0 & ...  & 0 \\
 1  & \frac{\binom{2}{1}\l^2}{1+\sqrt{1-\l^2}} & -\sqrt{1-\l^4} & ...  & 0 \\
 1  & \frac{\binom{3}{1}\l^2}{1+\sqrt{1-\l^2}} & \frac{\binom{3}{2}\l^4}{1+\sqrt{1-\l^4}} & ...  & 0 \\
 ... & ... & ... & ...  & ... \\
 1  & \frac{\binom{N}{1}\l^2}{1+\sqrt{1-\l^2}} & \frac{\binom{N}{2}\l^4}{1+\sqrt{1-\l^4}}  & ... & \frac{\binom{N}{N-1}\l^{2N-2}}{1+\sqrt{1-\l^{2N-2}}} \end{vmatrix}.
\]
One may also use the alternative recurrent formula
\[\lb{003b}
 \int_{-1}^1V(x)^Ndx=\frac1{\sqrt{1-\l^{2N}}}\lt(2+\sum_{n=1}^{N-1}\binom{N}{n}(1-\sqrt{1-\l^{2n}})\int_{-1}^1V(x)^ndx\rt).
\]
As an exercise, one can check
\[\lb{003cc}
 (1-\sqrt{1-\l^{2N}})\int_{-1}^{1}V(x)^Ndx=\int_{-1}^1(V(x)-1)^Ndx,\ \ \ N\ge1,
\]
and simplify \er{003} and \er{003a}, as something like
\[\lb{003ccc}
 \int_{-1}^{1}V(x)^Ndx=2\sum_{k\ge1}\sum_{N=N_0>...>N_k=0}(-1)^{N+k}\prod_{j=0}^{k-1}\binom{N_j}{N_{j+1}}\cdot\prod_{j=0}^{k-1}(1-\l^{2N_j})^{-\frac12}.
\] 
 
ii) Let $\cR_{\pm}h=h(\frac{x\pm1}2)$ be operators acting on $L^2(-1,1)$ (square integrable functions). Then
\[\lb{004}
 \int_{-1}^1 V(x)h(x) dx = \int_{-1}^1\oint_{|z|=1}\lt(1-\l\frac{z^{-1}\cR_-+z\cR_+}2\rt)^{-1} h(x)\frac{ dzdx}{2\pi i z}
\]
for any $h\in L^2$. Moreover, the $L^2\to L^2$-operator norm $\|z^{-1}\cR_-+z\cR_+\|\le2$ for $|z|=1$. Instead of $L^2$ one may take $L^{\iy}(-1,1)$ or $C([-1,1])$ (bounded or continuous functions).

iii) Let $P(x)=\sum_{n=0}^Np_nx^n$ be a polynomial with $p_n\in\C$ and even $N$. Then
\begin{multline}\lb{005}
	\int_{-1}^1 V(x)P(x) dx=\sum_{j=0}^{N}\frac{2}{2^{j}\sqrt{1-\frac{\l^2}{4^{j}}}\prod_{0\le n\ne j\le N}(1-\frac{2^n}{2^j})}\cdot\\
	{\begin{vmatrix} 1-\frac{2^0}{2^j} & 0 & 0 & 0 & ... & 0 & 1 \\ 
			-\binom{1}{0}\sqrt{1-\frac{\l^2}{4^{j}}} &  1-\frac{2^1}{2^j} & 0 &  0 & ... & 0 & 0 \\
			\binom{2}{0} & -\binom{2}{1}\sqrt{1-\frac{\l^2}{4^{j}}} &  1-\frac{2^2}{2^j} & 0 &... & 0 & \frac{2^2}3 \\
			-\binom{3}{0}\sqrt{1-\frac{\l^2}{4^{j}}} &  \binom{3}{1} & -\binom{3}{2}\sqrt{1-\frac{\l^2}{4^{j}}} &  1-\frac{2^3}{2^j} & ... & 0 & 0 \\
			... & ... & ... & ... & ... & ... & ... \\
			\binom{N}{0} & -\binom{N}{1}\sqrt{1-\frac{\l^2}{4^{j}}} & \binom{N}{2} & -\binom{N}{3}\sqrt{1-\frac{\l^2}{4^{j}}} & ... & 1-\frac{2^{N}}{2^j} & \frac{2^{N}}{N+1} \\
			-p_0 & -p_1 & -p_2 & -p_3 & ... & -p_N & 0
	\end{vmatrix}}.
\end{multline}
In particular, for even $N\ge0$ we have
\begin{multline}\lb{005a}
 \int_{-1}^1 V(x)x^N dx=\sum_{j=0}^{N}\frac{2}{2^{j}\sqrt{1-\frac{\l^2}{4^{j}}}\prod_{0\le n\ne j\le N}(1-\frac{2^n}{2^j})}\cdot\\
 {\begin{vmatrix} 1-\frac{2^0}{2^j} & 0 & 0 & 0 & ... & 1 \\ 
    	-\binom{1}{0}\sqrt{1-\frac{\l^2}{4^{j}}} &  1-\frac{2^1}{2^j} & 0 &  0 & ... & 0 \\
    	\binom{2}{0} & -\binom{2}{1}\sqrt{1-\frac{\l^2}{4^{j}}} &  1-\frac{2^2}{2^j} & 0 &... & \frac{2^2}3 \\
    	-\binom{3}{0}\sqrt{1-\frac{\l^2}{4^{j}}} &  \binom{3}{1} & -\binom{3}{2}\sqrt{1-\frac{\l^2}{4^{j}}} &  1-\frac{2^3}{2^j} & ... & 0 \\
    	... & ... & ... & ... & ... & ... \\
    	\binom{N}{0} & -\binom{N}{1}\sqrt{1-\frac{\l^2}{4^{j}}} & \binom{N}{2} & -\binom{N}{3}\sqrt{1-\frac{\l^2}{4^{j}}} & ... & \frac{2^{N}}{N+1}
    \end{vmatrix}}.
\end{multline}
If $N\ge0$ is odd then $\int_{-1}^1 V(x)x^N dx=0$.

iv) For $\o\in\C$ we have
\[\lb{006}
 \int_{-1}^1V(x)\cos\o xdx=\frac1{\pi}\int_{-\pi}^{\pi}C(\vp,\o)d\vp,
\]
where
\begin{multline}\lb{007}
 C(\vp,\o)=\sum_{n=0}^{\iy}\l^n\frac{2^n}{\o}\sin\frac{\o}{2^n}\prod_{j=1}^n\cos(\vp+\frac{\o}{2^j})=\\
 \frac{\sin\o}{\o}\cdot\cfrac{1}{1-\cfrac{\l\cos(\vp+\frac{\o}{2})}{\cos\frac{\o}{2}+\l\cos(\vp+\frac{\o}{2})-\cfrac{\l\cos\frac{\o}{2}\cos(\vp+\frac{\o}{4})}{\cos\frac{\o}{4}+\l\cos(\vp+\frac{\o}{4})-\cfrac{\l\cos\frac{\o}{4}\cos(\vp+\frac{\o}{8})}{\cos\frac{\o}{8}+\l\cos(\vp+\frac{\o}{8})-...}}}}.
\end{multline}
\end{theorem}
In Fig. \ref{fig2} we plot the Fourier series approximation of $V$, where the Fourier coefficients are computed by \er{006} and \er{007}.
Let us discuss the connection between ``fractal foothills'' $V$ and ``fractal mountains" $U$ defined in \cite{K}. Recall that $U$ is defined by
\[\lb{008}
 U(x)=\sum_{0\le m\le n<+\iy}\l^{n+1}L_{mn}(x),\ \ \ where\ \ \ L_{mn}(x)=\ca 1,& \sum_{j=m}^nx_j=0,\\ 0,& otherwise,\ac
\]
where $x_j$ are given in \er{001}. It is seen that $U$ counts all the loops in the ``random walk'' $x$, while $V$ counts the returns to the origin only, since $L_n(x)=L_{0n}(x)$, see \er{002} and \er{008}. It explains the fact why the structure of $V$ is much simpler than $U$. 
Using \er{001}, \er{002} and \er{008}, it is not difficult to write the explicit connection between $U$ and $V$, namely
\[\lb{009}
 U(x)=V(x)-1+\l(\cT V(x)-1)+\l^2(\cT^2V( x)-1)+...=(1-\l\cT)^{-1}V(x)-\frac1{1-\l},
\]
where $\cT$ is a change-of-variable operator that represents a left-shift of digits in the expansion \er{001}:
\[\lb{010}
 \cT V(x)=\ca V(2x-1), & x\in(0,1],\\ V(2x+1),& x\in[-1,0]. \ac
\]
Identity \er{009} is assumed to be valid in $L^2$, i.e., up to a set of zero Lebesgue measure. I made this remark to avoid the possible questions about including $x=0$ into the left or right interval. It is easy to check that adjoint operator $\cT^*=(\cR_++\cR_-)/2$, where $\cR_{\pm}$ are defined in Theorem \ref{T1}.ii. Using this fact along with \er{009} and the same ideas as in \er{123}-\er{125} for $z=1$ and for the basis $\{x^{2n}\}$ instead of $\{x^n\}$, we obtain statments i) and ii) of the following Corollary. Statement iii) is proven in the next Section.

\begin{corollary}\lb{C1}
For any $f\in L^2(-1,1)$, the following identity is fulfilled 
\[\lb{011}
 \int_{-1}^1U(x)f(x)dx=\int_{-1}^1V(x)\lt(1-\l\frac{\cR_++\cR_-}2\rt)^{-1}f(x)dx-\frac1{1-\l}\int_{-1}^1f(x)dx.
\]
For even polynomials $P(x)=\sum_{n=0}^{\frac{N}2}p_nx^{2n}$ with $p_n\in\C$ and even $N$, \er{011} implies
\begin{multline}\lb{012}
\int_{-1}^1U(x)P(x)dx=\frac{\begin{vmatrix} 1-\l & 0 & 0 & ... & 0 & \int_{-1}^1V(x)dx \\
		\frac{-\l\binom{2}{0}}{2^2} & 1-\frac{\l}{2^2} & 0 & ... & 0 & \int_{-1}^1V(x)x^2dx \\
		\frac{-\l\binom{4}{0}}{2^4} & \frac{-\l\binom{4}{2}}{2^4}  & 1-\frac{\l}{2^4} & ... & 0 &                      \int_{-1}^1V(x)x^4dx \\
		... & ... & ... & ... & ... & ... \\
		\frac{-\l\binom{N}{0}}{2^N} & \frac{-\l\binom{N}{2}}{2^N} & \frac{-\l\binom{N}{4}}{2^N} & ... & 1-\frac{\l}{2^N} &  \int_{-1}^1V(x)x^Ndx \\
		p_0 & p_1 & p_2 & ... & p_{\frac N2} & 0
\end{vmatrix}}{-\prod_{n=0}^{\frac{N}2}(1-\frac{\l}{4^n})}-\\
\sum_{n=0}^{\frac N2}\frac{2p_n}{(1-\l)(2n+1)},
\end{multline}
In particular, for even $N\ge0$ we have
\[\lb{012a}
 \int_{-1}^1U(x)x^Ndx=\frac{\begin{vmatrix} 1-\l & 0 & 0 & ... & \int_{-1}^1V(x)dx \\
                                      \frac{-\l\binom{2}{0}}{2^2} & 1-\frac{\l}{2^2} & 0 & ... & \int_{-1}^1V(x)x^2dx \\
 \frac{-\l\binom{4}{0}}{2^4} & \frac{-\l\binom{4}{2}}{2^4}  & 1-\frac{\l}{2^4} & ... &                      \int_{-1}^1V(x)x^4dx \\
 ... & ... & ... & ... & ... \\
 \frac{-\l\binom{N}{0}}{2^N} & \frac{-\l\binom{N}{2}}{2^N} & \frac{-\l\binom{N}{4}}{2^N} & ... &  \int_{-1}^1V(x)x^Ndx
                      \end{vmatrix}}{\prod_{n=0}^{\frac{N}2}(1-\frac{\l}{4^n})}-\frac{2}{(1-\l)(N+1)},
\]
where $\int_{-1}^1V(x)x^ndx$ can be computed by \er{005}. Note that if $N$ is odd then $\int_{-1}^1U(x)x^Ndx=0$, since $U$ is even function.

iii) The integration becomes simpler based on modified Bernoulli polynomials. Define $P_n(x):=2^nB_n(\frac{x+1}2)$, where $B_n$ are the classical Bernoulli polynomials. Then
\[\lb{013}
\int_{-1}^1U(x)P_n(x)dx=\frac{1}{1-2^{-n}\l}\int_{-1}^1V(x)P_n(x)dx-\frac{2\d_{n0}}{1-\l},\ \ n\ge0,
\]
where $\d$ is the Kronecker delta. In particular, for even $N$, we have
\begin{multline}\lb{014}
 \int_{-1}^1U(x)x^Ndx=\int_{-1}^1V(x)Q_N(x)dx-\frac{2}{(1-\l)(N+1)},\ \ where\\ Q_N(x)=\sum_{j=0}^{\frac N2}\frac{\binom{N}{2j}}{(1-2^{2j-N}\l)(2j+1)}P_{N-2j}(x).
\end{multline}
Denote $\cD=\frac{d}{dx}$. There are a few useful relations for the polynomials $P_n(x)$:
\[\lb{015}
 \sum_{n=0}^{+\iy}P_n(x)\frac{t^n}{n!}=\frac{t}{\sinh t}e^{tx},\ \ \ P_n(x)=\frac{\cD}{\sinh\cD}x^n,\ \ \ \cD P_n(x)=nP_{n-1}(x).
\]
\end{corollary} 

{\bf Remark.} Polynomials $\{P_n\}_{n\ge0}$ is an Appell sequence, since $\cD P_n(x)=nP_{n-1}(x)$, see \er{015}. Formula $P_n(x)=\frac{\cD}{\sinh\cD}x^n$ is convenient for calculating $P_n(x)$. We have
\[\lb{016}
 P_n(x)=\sum_{j=0}^n\binom{n}{j}c_jx^{n-j},\ \ with\ \ \frac{\cD}{\sinh \cD}=\sum_{n=0}^{+\iy}\frac{c_n}{n!}\cD^n.
\]
Thus, all $c_{2n+1}=0$ and
\[\lb{017}
 c_0=1,\ c_2=\frac{-1}{3},\ c_4=\frac7{15},\ c_6=\frac{-31}{21},\ c_8=\frac{127}{15},\ ...,\ c_{2n}=-\sum_{j=0}^{n-1}\frac{\binom{2n}{2j}c_{2j}}{2n-2j+1}.
\]
Further analysis may be based on \er{006}, \er{007}, and new formula
\[\lb{018}
 \int_{-1}^1U(x)(e^{\o x}-\l\cosh\frac{\o}2e^{\frac{\o x}2})dx=\int_{-1}^1V(x)e^{\o x}dx-\frac{2\sinh\o}{\o},\ \ \o\in\C
\]
that immediately follows from \er{011}.

We have obtained \er{012} as the alternative formula to the already presented one in \cite{K}.  At the same time, the closed form expression for $\int_{-1}^1 U(x)^Ndx$ similar to \er{003}, \er{003a} and \er{003b} is still a good challenge, at least to me. I believe also that there are further simplifications of \er{005} and \er{012}, not obvious to me at the moment.

\begin{figure}
  \centering
    \includegraphics[width=0.99\textwidth]{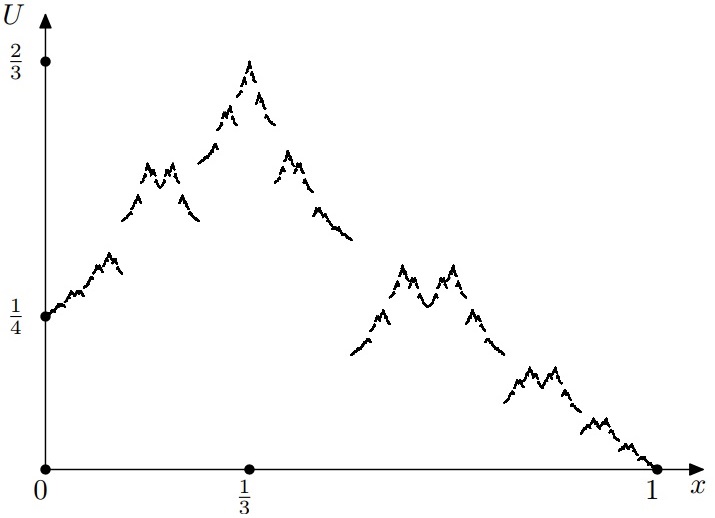}
  \caption{The plot of $U$ for $\l=1/2$. \lb{fig3}}  
\end{figure}
\begin{figure}
    \centering
    \begin{subfigure}[b]{0.49\textwidth}
        \includegraphics[width=\textwidth]{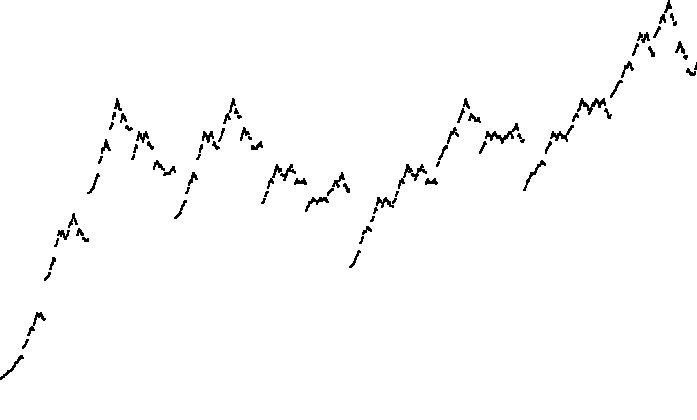}
    \end{subfigure}
    \begin{subfigure}[b]{0.49\textwidth}
        \includegraphics[width=\textwidth]{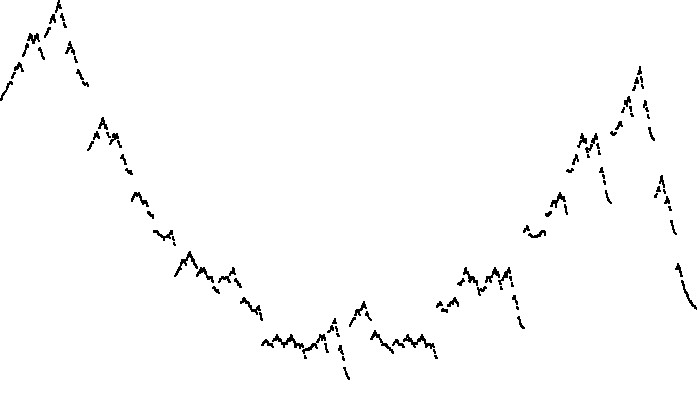}
    \end{subfigure}
    \caption{For $\l=1/2$, two randomly chosen different segments of the curve $U$ zoomed in $\approx2^{70}$ times. 
}\label{fig4}
\end{figure}

Let us provide a few formulas followed from Theorem \ref{T1} and Corollary \ref{C1}. This also reduces some disambiguation in reading \er{005}, \er{005a} and \er{012}, \er{012a} for small $N$ ($N=0$). We have
$$
 \int_{-1}^1 V(x)dx=\frac2{\sqrt{1-\l^2}},\ \ \ \int_{-1}^1V(x)^2dx=\frac4{\sqrt{1-\l^4}\sqrt{1-\l^2}}-\frac2{\sqrt{1-\l^4}},
$$
$$
 \int_{-1}^1V(x)x^2dx=\frac{4\sqrt{1-\l^2}}{3}+\frac2{3\sqrt{1-\l^2}}-4\sqrt{1-\frac{\l^2}4}+\frac{8\sqrt{1-\frac{\l^2}{16}}}3,
$$
$$
 \int_{-1}^1U(x)dx=\frac{2}{(1-\l)\sqrt{1-\l^2}}-\frac2{1-\l},
$$
$$
 \int_{-1}^1U(x)x^2=\frac{\frac{4\sqrt{1-\l^2}}3-4\sqrt{1-\frac{\l^2}4}+\frac{8\sqrt{1-\frac{\l^2}{16}}}3}{1-\frac{\l}{4}}+\frac{2}{3(1-\l)\sqrt{1-\l^2}}-\frac{2}{3(1-\l)},
$$
where the last two integrals are already presented in \cite{K}. Let us conclude with a few words about the comparison of $V(x)$ and $U(x)$. The first function is already a fractal curve, but the second one is a ``double" fractal curve, since we apply the ``fractal" resize-operator $\cT$ to the already fractal curve $V$, see \er{009} and \er{010}. We can compare the plots of $V$, see Fig. \ref{fig1}, and $U$ presented on Figs. \ref{fig3} and \ref{fig4}. The first plot I have taken from \cite{K}, but the zoomed ones are new.

{\section{Proof of the main results\lb{sec2}}

\subsection{Analytic generating function for $V$.} For $x\in[-1,1]$, let us define the function
\[\lb{101}
 F(x,z)=1+\l z^{x_0}+\l^2z^{x_0+x_1}+\l^3z^{x_0+x_1+x_2}+...=1+\sum_{n=0}^{+\iy}\l^{n+1}z^{\sum_{j=0}^nx_j},
\]
where $x_n\in\{-1,+1\}$ are given by \er{001}. Since $|\l|<1$, it is seen that for any $x\in[-1,1]$ function $F$ is analytic in some open ring containing the circle $|z|=1$. Indeed, each term of the series \er{101} can be uniformly approximated by the terms of a convergent series
\[\lb{102}
 |F(x,z)|\le 1+\l\max\{|z|,|z|^{-1}\}+\l^2\max\{|z|^2,|z|^{-2}\}+...\le\frac1{1-\l|z|}+\frac1{1-\l|z|^{-1}},
\]
since all $x_n\in\{-1,+1\}$. Thus $F(x,z)$ is analytic in $\{z:\ \l<|z|<\l^{-1}\}$ for any fixed $x\in[-1,1]$. It is seen that $V$ is a free term in the ($z$-)series for $F$, see \er{002} and \er{101}. Thus, we have
\[\lb{103}
 V(x)=\mathfrak{f}(F(x,z))=\oint_{|z|=1}F(x,z)\frac{dz}{2\pi i z},
\]
where symbol $\mathfrak{f}$ means the free term in the Laurent series. Using \er{001} and \er{101}, we derive the functional equation
\[\lb{104}
 F(\pm\frac12+y,z)=1+\l z^{\pm1}F(2y,z),\ \ \ y\in[-\frac12,\frac12],
\]
basic in our research. 

\subsection{Integrals $\int V(x)^Ndx$.} There are many possible ways, we chose an exotic one. For $\a\ss\{1,...,N\}$, let us denote
\[\lb{105}
 F_{\a}:= \prod_{j\in\a}F(x,z_j),\ \ \ z_{\a}:=\prod_{j\in\a} z_{j},\ \ \ F_{\es}=z_{\es}=1.
\]
Using \er{104}, we obtain
\[\lb{106}
 \cJ F_{\a}=\frac12\cJ\prod_{j\in\a}(1+\l z_jF_j)+\frac12\cJ\prod_{j\in\a}(1+\l z_j^{-1}F_j)=\sum_{\b\ss\a}\l^{|\b|}\frac{z_{\b}+z_{\b}^{-1}}2\cJ F_{\b}
\]
and, hence,
\[\lb{107}
 \cJ F_{\a}=\frac1{1-\l^{|\a|}\frac{z_{\a}+z_{\a}^{-1}}2}\sum_{\b\subsetneqq\a}\l^{|\b|}\frac{z_{\b}+z_{\b}^{-1}}2\cJ F_{\b},
\]
where here and below $\cJ\cdot:=\int_{-1}^{1}\cdot dx$ and $|\cdot|$ denotes the number of elements. Thus, we deduce that
\[\lb{108}
 \cJ F_{\a}=\frac1{1-\l^{|\a|}\frac{z_{\a}+z_{\a}^{-1}}2}\lt(\cJ1+\sum_{\es\subsetneqq\g\subsetneqq...\subsetneqq\b\subsetneqq\a}\l^{|\b|}\frac{\frac{z_{\b}+z_{\b}^{-1}}2}{1-\l^{|\b|}\frac{z_{\b}+z_{\b}^{-1}}2}...\l^{|\g|}\frac{\frac{z_{\g}+z_{\g}^{-1}}2}{1-\l^{|\g|}\frac{z_{\g}+z_{\g}^{-1}}2}\cJ1\rt),
\]
where the sum is taken over all nested sequences of sets $\es\subsetneqq\g\subsetneqq...\subsetneqq\b\subsetneqq\a$, and $\cJ1=2$. Since the sets are embedded strictly in each other, we conclude that the free term of the Laurent ($z_1,...,z_N$-)series \er{108} satisfies the equality
\[\lb{109}
 \mathfrak{f}(\cJ F_{\a})=\mathfrak{f}\lt(\frac1{1-\l^{|\a|}\frac{z_{\a}+z_{\a}^{-1}}2}\rt)\mathfrak{f}\lt(2+2\sum_{\es\subsetneqq\g\subsetneqq...\subsetneqq\b\subsetneqq\a}\l^{|\b|}\frac{\frac{z_{\b}+z_{\b}^{-1}}2}{1-\l^{|\b|}\frac{z_{\b}+z_{\b}^{-1}}2}...\l^{|\g|}\frac{\frac{z_{\g}+z_{\g}^{-1}}2}{1-\l^{|\g|}\frac{z_{\g}+z_{\g}^{-1}}2}\rt)
\]
and, similarly,
\begin{multline}\lb{110}
\mathfrak{f}\lt(\l^{|\a|}{\frac{z_{\a}+z_{\a}^{-1}}2}\cJ F_{\a}\rt)=\mathfrak{f}\lt(\frac{\l^{|\a|}{\frac{z_{\a}+z_{\a}^{-1}}2}}{1-\l^{|\a|}\frac{z_{\a}+z_{\a}^{-1}}2}\rt)\cdot\\
\mathfrak{f}\lt(2+2\sum_{\es\subsetneqq\g\subsetneqq...\subsetneqq\b\subsetneqq\a}\l^{|\b|}\frac{\frac{z_{\b}+z_{\b}^{-1}}2}{1-\l^{|\b|}\frac{z_{\b}+z_{\b}^{-1}}2}...\l^{|\g|}\frac{\frac{z_{\g}+z_{\g}^{-1}}2}{1-\l^{|\g|}\frac{z_{\g}+z_{\g}^{-1}}2}\rt).
\end{multline}
Using the facts that
\begin{multline}\lb{111}
 \mathfrak{f}\lt(\frac1{1-\l^{|\a|}\frac{z_{\a}+z_{\a}^{-1}}2}\rt)=\mathfrak{f}\lt(\frac1{1-\l^{|\a|}\frac{z+z^{-1}}2}\rt)=\oint_{|z|=1}\frac1{1-\l^{|\a|}\frac{z+z^{-1}}2}\frac{dz}{2\pi i z}=\\
 {\rm Res}\lt(\frac2{2z-\l^{|\a|}z^2-\l^{|\a|}};z=\frac{1-\sqrt{1-\l^{2|\a|}}}{\l^{|\a|}}\rt)=\frac1{\sqrt{1-\l^{2|\a|}}}
\end{multline}
and
\[\lb{112}
 \mathfrak{f}\lt(\frac{\l^{|\a|}{\frac{z_{\a}+z_{\a}^{-1}}2}}{1-\l^{|\a|}\frac{z_{\a}+z_{\a}^{-1}}2}\rt)=\mathfrak{f}\lt(-1+\frac1{1-\l^{|\a|}\frac{z_{\a}+z_{\a}^{-1}}2}\rt)=\frac1{\sqrt{1-\l^{2|\a|}}}-1
\]
along with \er{109} and\er{110}, we deduce that
\[\lb{113}
 \mathfrak{f}\lt(\l^{|\a|}{\frac{z_{\a}+z_{\a}^{-1}}2}\cJ F_{\a}\rt)=(1-\sqrt{1-\l^{2|\a|}})\mathfrak{f}(\cJ F_{\a}).
\]
Thus, using \er{113}, \er{107}, and simple fact that $\mathfrak{f}(\cJ F_{\a})=\mathfrak{f}(\cJ F_{\a'})$ if $|\a|=|\a'|$, we obtain
\[\lb{114}
 \mathfrak{f}(\cJ F_{\{1,...,N\}})=\frac1{\sqrt{1-\l^{2N}}}\lt(2+\sum_{n=1}^{N-1}\binom{N}{n}(1-\sqrt{1-\l^{2n}})\mathfrak{f}(\cJ F_{\{1,...,n\}})\rt)
\]
that can be written in the matrix form
$$
 \ma
   1 & 0 & 0 & ... & 0 \\
   -1  & \sqrt{1-\l^2} & 0 & ... & 0 \\
   -1  & \binom{2}{1}(\sqrt{1-\l^2}-1) & \sqrt{1-\l^4} & ... & 0 \\
   ... & ... & ... & ... & ... \\
   -1  & \binom{N}{1}(\sqrt{1-\l^2}-1) & \binom{N}{2}(\sqrt{1-\l^4}-1)  & ... & \sqrt{1-\l^{2N}}  
 \am
 \ma \mathfrak{f}(\cJ F_{\es}) \\ \mathfrak{f}(\cJ F_{\{1\}}) \\ \mathfrak{f}(\cJ F_{\{1,2\}}) \\ ... \\ \mathfrak{f}(\cJ F_{\{1,...,N\}}) \am =\ma 2 \\ 0 \\ 0 \\ ... \\ 0 \am
$$
or
\[\lb{115}
 \ma
   1 & 0 & 0 & ... & 0 \\
   1  & -\sqrt{1-\l^2} & 0 & ... & 0 \\
   1  & \frac{\binom{2}{1}\l^2}{1+\sqrt{1-\l^2}} & -\sqrt{1-\l^4} & ... & 0 \\
   ... & ... & ... & ... & ... \\
   1  & \frac{\binom{N}{1}\l^2}{1+\sqrt{1-\l^2}} & \frac{\binom{N}{2}\l^2}{1+\sqrt{1-\l^2}}  & ... & -\sqrt{1-\l^{2N}}  
 \am
 \ma \mathfrak{f}(\cJ F_{\es}) \\ \mathfrak{f}(\cJ F_{\{1\}}) \\ \mathfrak{f}(\cJ F_{\{1,2\}}) \\ ... \\ \mathfrak{f}(\cJ F_{\{1,...,N\}}) \am =\ma 2 \\ 0 \\ 0 \\ ... \\ 0 \am.
\]
Applying the general Cramer's rule to the linear system \er{115} we get the identity
\begin{multline}\lb{116}
p_0\mathfrak{f}(\cJ F_{\es})+...+p_N\mathfrak{f}(\cJ F_{\{1,...,N\}})=\\
\frac{-(-1)^N}{\prod_{n=1}^N\sqrt{1-\l^{2n}}}\begin{vmatrix}
   1 & 0 & 0 & ... & 0 & 2 \\
   1  & -\sqrt{1-\l^2} & 0 & ... & 0 & 0 \\
   1  & \frac{\binom{2}{1}\l^2}{1+\sqrt{1-\l^2}} & -\sqrt{1-\l^4} & ... & 0 & 0 \\
   ... & ... & ... & ... & ... \\
   1  & \frac{\binom{N}{1}\l^2}{1+\sqrt{1-\l^2}} & \frac{\binom{N}{2}\l^4}{1+\sqrt{1-\l^4}}  & ... & -\sqrt{1-\l^{2N}} & 0 \\
   p_0 & p_1 & p_2 & ... & p_N & 0  
 \end{vmatrix},\ \ \ p_i\in\C.
\end{multline}
Finally, expanding \er{116} by the last column and using the obvious extension of \er{103} based on \er{105}, namely,
\[\lb{117}
 \int_{-1}^{1}V(x)^Ndx=\oint_{|z_1|=1}...\oint_{|z_N|=1}\int_{-1}^1F_{\{1,...,N\}}dx\frac{dz_N}{2\pi i z_N}...\frac{dz_1}{2\pi i z_1}=\mathfrak{f}(\cJ F_{\{1,...,N\}}),
\]
we obtain the announced formula \er{003}. Formula \er{003a} is a simple consequence of \er{003}. The alternative recurrent formula \er{003b} follows from \er{114} and \er{117}. 

\subsection{Operator's identity for $\int V(x)h(x)dx$, $h\in L^2$.} We use the notation $\cJ:=\int_{-1}^1\cdot dx$ from the previous subsection. Using \er{104} by analogy with \er{106}, we obtain
\[\lb{118}
 \cJ F\wt h=\cJ \wt h+\l\cJ F\frac{z^{-1}\cR_-+z\cR_+}2\wt h,\ \ \forall\wt h\in L^2(-1,1),
\]
that leads to
\[\lb{119}
 \cJ F\lt(1-\l\frac{z^{-1}\cR_-+z\cR_+}2\rt)\wt h=\cJ \wt h
\]
or
\[\lb{120}
 \cJ F h=\cJ \lt(1-\l\frac{z^{-1}\cR_-+z\cR_+}2\rt)^{-1} h,\ \ \forall h\in L^2(-1,1)
\]
if the corresponding inverse operator in \er{120} exists. To show the existence of the inverse operator in \er{120}, it is enough to show that the operator norm $\|z^{-1}\cR_-+z\cR_+\|_{L^2\to L^2}\le 2$, because using $|\l|<1$ we may write the converging geometric series for the inverse operator. The mentioned norm's inequality follows from
\begin{multline}\lb{121}
 \|zf(\frac{x+1}2)+z^{-1}f(\frac{x-1}2)\|_{L^2(-1,1)}^2\le\|zf(\frac{x+1}2)+
 z^{-1}f(\frac{x-1}2)\|_{L^2(-1,1)}^2+\\
 \|zf(\frac{x+1}2)-z^{-1}f(\frac{x-1}2)\|_{L^2(-1,1)}^2=2\|f(\frac{x+1}2)\|_{L^2(-1,1)}^2+2\|f(\frac{x-1}2)\|_{L^2(-1,1)}^2=\\
 4\|f(x)\|^2_{L^2(0,1)}+4\|f(x)\|^2_{L^2(-1,0)}=4\|f(x)\|^2_{L^2(-1,1)},
\end{multline}
where we use the fact that $|z|=1$. Also, it is much easy to check that the operator norm $\|z^{-1}\cR_-+z\cR_+\|_{L^{\iy}\to L^{\iy}}\le 2$ if we consider the space of bounded functions $L^{\iy}$ (or $C$) instead of $L^2$. Now, the announced formula \er{004} follows from \er{120} and 
\[\lb{122}
 \int_{-1}^1 V(x)h(x)dx=\oint_{|z|=1}\cJ Fh\frac{dz}{2\pi i z},
\]
that is similar to \er{103} and \er{117}.

\subsection{Integrals $\int V(x)x^Ndx$.} At first, let us compute
\begin{multline}\lb{123}
 \lt(1-\l\frac{z^{-1}\cR_-+z\cR_+}2\rt)x^n=x^n-\frac{\l z}2\lt(\frac{x+1}2\rt)^n-\frac{\l z^{-1}}2\lt(\frac{x-1}2\rt)^n=\\
 (1-\frac{\l(z+z^{-1})}{2^{n+1}})x^n-\frac{\l}{2^{n+1}}\sum_{j=0}^{n-1} \binom{n}{j}((-1)^{n-j}z^{-1}+z)x^{j},
\end{multline}
that can be written in the matrix form
\begin{multline}\lb{124}
\lt(1-\l\frac{z^{-1}\cR_-+z\cR_+}2\rt)\ma 1 \\ x \\ x^2 \\ ... \\ x^N \am= \\
\ma 1-\frac{\l(z+z^{-1})}{2} & 0 & 0 & ... & 0 \\ 
    \frac{-\l\binom{1}{0}((-1)^{1-0}z^{-1}+z)}{2^{2}} & 1-\frac{\l(z+z^{-1})}{2^2} & 0 & ... & 0 \\
    \frac{-\l\binom{2}{0}((-1)^{2-0}z^{-1}+z)}{2^{3}} & \frac{-\l\binom{2}{1}((-1)^{2-1}z^{-1}+z)}{2^{3}} & 1-\frac{\l(z+z^{-1})}{2^3} & ... & 0 \\
    ... & ... & ... & ... & ... \\
    \frac{-\l\binom{N}{0}((-1)^{N-0}z^{-1}+z)}{2^{N+1}} & \frac{-\l\binom{N}{1}((-1)^{N-1}z^{-1}+z)}{2^{N+1}} & \frac{-\l\binom{N}{2}((-1)^{N-2}z^{-1}+z)}{2^{N+1}} & ... & 1-\frac{\l(z+z^{-1})}{2^{N+1}}
\am
\ma 1 \\ x \\ x^2 \\ ... \\ x^N \am.
\end{multline}
Thus, using Cramer's rule, we write
\begin{multline}\lb{125}
\lt(1-\l\frac{z^{-1}\cR_-+z\cR_+}2\rt)^{-1}\sum_{n=0}^Np_nx^n=\ma p_0 & p_1 & p_2 & ... & p_N \am\lt(1-\l\frac{z^{-1}\cR_-+z\cR_+}2\rt)^{-1}\ma 1 \\ x \\ x^2 \\ ... \\ x^N \am= \\
\frac{-1}{\prod_{n=1}^{N+1}(1-\frac{\l(z+z^{-1})}{2^n})}\cdot
\\
\begin{vmatrix} 1-\frac{\l(z+z^{-1})}{2} & 0 & 0 & ... & 0 & 1 \\ 
    \frac{-\l\binom{1}{0}((-1)^{1-0}z^{-1}+z)}{2^{2}} & 1-\frac{\l(z+z^{-1})}{2^2} & 0 & ... & 0 & x \\
    \frac{-\l\binom{2}{0}((-1)^{2-0}z^{-1}+z)}{2^{3}} & \frac{-\l\binom{2}{1}((-1)^{2-1}z^{-1}+z)}{2^{3}} & 1-\frac{\l(z+z^{-1})}{2^3} & ... & 0 & x^2 \\
    ... & ... & ... & ... & ... & ... \\
    \frac{-\l\binom{N}{0}((-1)^{N-0}z^{-1}+z)}{2^{N+1}} & \frac{-\l\binom{N}{1}((-1)^{N-1}z^{-1}+z)}{2^{N+1}} & \frac{-\l\binom{N}{2}((-1)^{N-2}z^{-1}+z)}{2^{N+1}} & ... & 1-\frac{\l(z+z^{-1})}{2^{N+1}} & x^N \\
    p_0 & p_1 & p_2 & ... & p_N & 0
\end{vmatrix}=\\
\frac{-2}{\prod_{n=1}^{N+1}(\l z^2-2^nz+\l)}\cdot
\\
\begin{vmatrix} \l z^2-2z+\l & 0 & 0 & ... & 0 & 1 \\ 
    \l\binom{1}{0}((-1)^{1-0}+z^2) &  \l z^2-2^2z+\l & 0 & ... & 0 & 2x \\
    \l\binom{2}{0}((-1)^{2-0}+z^2) & \l\binom{2}{1}((-1)^{2-1}+z^2) &  \l z^2-2^3z+\l & ... & 0 & (2x)^2 \\
    ... & ... & ... & ... & ... & ... \\
    \l\binom{N}{0}((-1)^{N-0}+z^2) & \l\binom{N}{1}((-1)^{N-1}+z^2) & \l\binom{N}{2}((-1)^{N-2}+z^2) & ... & \l z^2-2^{N+1}z+\l & (2x)^N \\
    -zp_0 & -zp_1 & -zp_2 & ... & -zp_N & 0
\end{vmatrix}.
\end{multline}
Using \er{125}, we write
\begin{multline}\lb{126}
	\cJ\frac1z\lt(1-\l\frac{z^{-1}\cR_-+z\cR_+}2\rt)^{-1}\sum_{n=0}^Np_nx^n=
	\frac{-2}{\prod_{n=1}^{N+1}(\l z^2-2^{n}z+\l)}\cdot
	\\ {\footnotesize
	\begin{vmatrix} \l z^2-2z+\l & 0 & 0 & ... & 0 & 2 \\ 
		\l\binom{1}{0}(z^2-1) &  \l z^2-2^2z+\l & 0 & ... & 0 & 0 \\
		\l\binom{2}{0}(z^2+1) & \l\binom{2}{1}(z^2-1) &  \l z^2-2^3z+\l & ... & 0 & \frac{2^3}3 \\
		... & ... & ... & ... & ... & ... \\
		\l\binom{N}{0}((-1)^{N-0}+z^2) & \l\binom{N}{1}((-1)^{N-1}+z^2) & \l\binom{N}{2}((-1)^{N-2}+z^2) & ... & \l z^2-2^{N+1}z+\l & \frac{2^N(1+(-1)^{N})}{N+1} \\
		-p_0 & -p_1 & -p_2 & ... & -p_N & 0
	\end{vmatrix}}.
\end{multline}
Now, we need to compute the integral of \er{126} over the unit circle. All poles of RHS in \er{126} that lie inside the unit ball are simple, they are smaller than the roots of polynomials $2^{n}z-\l z^2-\l$ and have the form
\[\lb{127}
 z_n=\frac{2^{n-1}-\sqrt{4^{n-1}-\l^2}}{\l},\ \ \ n=1,...,N+1.
\]
We assume that $N$ is even.  Applying the Cauchy residue theorem to \er{126} and using identities
\[\lb{126a}
 \l(z_n^2+1)=2^nz_n,\ \ \ \l(z_n^2-1)=-2\sqrt{4^{n-1}-\l^2}z_n,
\]
we deduce that
\begin{multline}\lb{127a}
	\cJ\oint_{|z|=1}\lt(1-\l\frac{z^{-1}\cR_-+z\cR_+}2\rt)^{-1}\sum_{n=0}^Np_nx^n\frac{dz}{2\pi iz}=
	\sum_{j=1}^{N+1}\frac{1}{\sqrt{4^{j-1}-\l^2}\prod_{1\le n\ne j\le N+1}z_j(2^{j}-2^n)}\cdot
	\\
	\begin{vmatrix} z_j(2^j-2) & 0 & 0 & ... & 0 & 2 \\ 
		-2z_j\binom{1}{0}\sqrt{4^{j-1}-\l^2} &  z_j(2^j-2^2) & 0 & ... & 0 & 0 \\
		2^jz_j\binom{2}{0} & -2z_j\binom{2}{1}\sqrt{4^{j-1}-\l^2} &  z_j(2^j-2^3) & ... & 0 & \frac{2^3}3 \\
		... & ... & ... & ... & ... & ... \\
		2^jz_j\binom{N}{0} & -2z_j\binom{N}{1}\sqrt{4^{j-1}-\l^2} & 2^jz_j\binom{N}{2} & ... & z_j(2^j-2^{N+1}) & \frac{2^{N+1}}{N+1} \\
		-p_0 & -p_1 & -p_2 & ... & -p_N & 0
	\end{vmatrix}=\\
    \sum_{j=1}^{N+1}\frac{1}{\sqrt{4^{j-1}-\l^2}\prod_{1\le n\ne j\le N+1}(1-2^{n-j})}\cdot
    \\
    \begin{vmatrix} 1-2^{1-j} & 0 & 0 & ... & 0 & 2 \\ 
    	-\binom{1}{0}\sqrt{1-\frac{\l^2}{4^{j-1}}} &  1-2^{2-j} & 0 & ... & 0 & 0 \\
    	\binom{2}{0} & -\binom{2}{1}\sqrt{1-\frac{\l^2}{4^{j-1}}} &  1-2^{3-j} & ... & 0 & \frac{2^3}3 \\
    	... & ... & ... & ... & ... & ... \\
    	\binom{N}{0} & -\binom{N}{1}\sqrt{1-\frac{\l^2}{4^{j-1}}} & \binom{N}{2} & ... & 1-2^{N+1-j} & \frac{2^{N+1}}{N+1} \\
    	-p_0 & -p_1 & -p_2 & ... & -p_N & 0
    \end{vmatrix},
\end{multline}
which with \er{004} give \er{005}. Formula \er{005a} follows from \er{005} directly.

\subsection{Fourier transform of $V$.} Again, we apply \er{004}. Using the fact that $|\l|<1$, and $\|z^{-1}\cR_-+z\cR_+\|\le2$ for $|z|=1$, we obtain the following geometric series expansion
\begin{multline}\lb{128}
 \lt(1-\l\frac{e^{-i\vp}\cR_-+e^{i\vp}\cR_+}2\rt)^{-1}e^{i\o x}=\sum_{n=0}^{\iy}\l^n\lt(\frac{e^{-i\vp}\cR_-+e^{i\vp}\cR_+}2\rt)^ne^{i\o x}=\\
\sum_{n=0}^{\iy}\l^n\lt(\frac{e^{-i\vp-i\frac{\o}2}+e^{i\vp+i\frac{\o}2}}2\rt)...\lt(\frac{e^{-i\vp-i\frac{\o}{2^n}}+e^{i\vp+i\frac{\o}{2^n}}}2\rt)e^{i\frac{\o x}{2^n}}= \sum_{n=0}^{\iy}e^{i\frac{\o x}{2^n}}\l^n\prod_{j=1}^n\cos(\vp+\frac{\o}{2^j}),
\end{multline}
where $\vp\in\R$ and $\o\in\C$. Thus, \er{128} along with \er{004} give us
\[\lb{129}
 \int_{-1}^1e^{i\o x}V(x)dx=\cJ\int_{-\pi}^{\pi}\lt(1-\l\frac{e^{-i\vp}\cR_-+e^{i\vp}\cR_+}2\rt)^{-1}e^{i\o x}\frac{d\vp}{2\pi}=\frac1{\pi}\int_{-\pi}^{\pi}C(\vp,\o),
\] 
where
\begin{multline}\lb{130}
 C(\vp,\o)=\sum_{n=0}^{\iy}\l^n\frac{2^n}{\o}\sin\frac{\o}{2^n}\prod_{j=1}^n\cos(\vp+\frac{\o}{2^j})=\\
 \frac{\sin\o}{\o}\sum_{n=0}^{\iy}\prod_{j=1}^n\frac{\l\cos(\vp+\frac{\o}{2^j})}{\cos\frac{\o}{2^j}}=\frac{\sin\o}{\o}\cdot\cfrac{1}{1-\cfrac{\frac{\l\cos(\vp+\frac{\o}{2})}{\cos\frac{\o}{2}}}{1+\frac{\l\cos(\vp+\frac{\o}{2})}{\cos\frac{\o}{2}}-\cfrac{\frac{\l\cos(\vp+\frac{\o}{4})}{\cos\frac{\o}{4}}}{ 1+\frac{\l\cos(\vp+\frac{\o}{4})}{\cos\frac{\o}{4}}-\cfrac{\frac{\l\cos(\vp+\frac{\o}{8})}{\cos\frac{\o}{8}}}{1+\frac{\l\cos(\vp+\frac{\o}{8})}{\cos\frac{\o}{8}}-...}}}}=\\
 \frac{\sin\o}{\o}\cdot\cfrac{1}{1-\cfrac{\l\cos(\vp+\frac{\o}{2})}{\cos\frac{\o}{2}+\l\cos(\vp+\frac{\o}{2})-\cfrac{\l\cos\frac{\o}{2}\cos(\vp+\frac{\o}{4})}{\cos\frac{\o}{4}+\l\cos(\vp+\frac{\o}{4})-\cfrac{\l\cos\frac{\o}{4}\cos(\vp+\frac{\o}{8})}{\cos\frac{\o}{8}+\l\cos(\vp+\frac{\o}{8})-...}}}}.
\end{multline}
To derive last three identities in \er{130}, we have used the identity $\sin\o=2^n\sin\frac{\o}{2^n}\prod_{j=1}^m\cos\frac{\o}{2^j}$ and Euler's continued fraction formula. Note that while the first identity in \er{130} is valid for $\vp\in\R$, $\o\in\C$ (for $\o=0$ there is a limit $\frac{\sin\o}{\o}\to1$), the other identities in \er{130} are formally valid for $\vp\in\R$, $\o\in\C\sm\pi\Z$. Using \er{129} and \er{130} and the fact that $V$ is an even function, we obtain the announced formulas \er{006} and \er{007}.

\subsection{Proof of Corollary \ref{C1}.iii).} Using the generating function for Bernoulli polynomials, we obtain
\[\lb{131}
 P(t,x):=\sum_{n=0}^{+\iy}P_n(x)\frac{t^n}{n!}=\sum_{n=0}^{+\iy}2^nB_n(\frac{x+1}2)\frac{t^n}{n!}=\frac{2te^{2t\frac{x+1}2}}{e^{2t}-1}=\frac{t}{\sinh t}e^{tx},
\]
which with
\[\lb{132}
 (\cR_-+\cR_+)P(t,x)=\frac{t}{\sinh t}(\cR_-+\cR_+)e^{tx}=\frac{t}{\sinh t}(e^{\frac t2}+e^{-\frac t2})e^{\frac{tx}2}=2P(\frac t2,x)
\]
leads to
\[\lb{133}
 (\cR_-+\cR_+)P_n(x)=\frac{2}{2^n}P_n(x).
\]
Using identity $\cD^ke^{tx}=t^ke^{tx}$, $k\ge0$ along with \er{131}, we obtain
\[\lb{134}
 P(t,x)=\frac{\cD}{\sinh\cD}e^{tx}.
\]
Expanding $e^{tx}$ into the Taylor series and using \er{134} along with the definition \er{131} we obtain
\[\lb{135}
 P_n(x)=\frac{\cD}{\sinh\cD}x^n
\]
that also leads to
\[\lb{136}
 \cD P_n(x)=\cD\frac{\cD}{\sinh\cD}x^n=\frac{\cD}{\sinh\cD}\cD x^n=n\frac{\cD}{\sinh\cD}x^{n-1}=n P_{n-1}(x).
\]
Both identities \er{135} and \er{136} for even $n$ give us
\[\lb{137}
 x^n=\frac{\sinh \cD}{\cD}P_n(x)=\sum_{j=0}^{\frac n2}\frac1{(2j+1)!}\cD^{2j}P_{n}(x)=\sum_{j=0}^{\frac n2}\frac{\binom{n}{2j}}{2j+1}P_{n-2j}(x).
\] 
Now, all the ingredients are ready. Formula \er{013} follows from \er{011} and \er{133}. Formula \er{014} follows from \er{013} and \er{137}. Formulas \er{131}, \er{135} and \er{136} imply \er{015}.

\section*{Conflicts of interests} There are no conflicts of interests or competing interests known to the author.

\section*{Data availability statement} The data that support the findings of this study are available from the corresponding author, upon reasonable request.

\end{document}